\documentclass[12pt,letterpaper]{article}
\usepackage[latin1]{inputenc}
\usepackage{amsmath}
\usepackage{amsfonts}
\usepackage{amssymb}
\usepackage{graphicx}
\usepackage{cite}
\usepackage[title]{appendix}
\author{Kyle Johnson\\
	\textit{Department of Mathematics}\\
	\textit{University of California}\\
	\textit{Davis, CA 95616, USA}\\
	\textit{email: kyljohnson@ucdavis.edu}}

\newcommand{\x}{\xi}
\newcommand{\la}{\lambda}
\newcommand{\e}{\epsilon}
\newcommand{\g}{\gamma}
\newcommand{\s}{\sigma}
\newcommand{\z}{{\zeta}}
\newcommand{\ph}{\varphi}
\newcommand{\G}{\Gamma}
\newcommand{\cC}{\mathcal{C}}
\newcommand{\bbP}{\mathbb{P}}
\newcommand{\bbC}{\mathbb{C}}
\newcommand{\Z}{\mathbb{Z}}

\newcommand{\tr}{\text{tr}}
\newcommand{\Ai}{\text{Ai}}
\newcommand{\Airy}{\text{Airy}}
\newcommand{\wsum}{\sum_{j=1}^L\frac{w_j}{w_j\xi-\tau}}
\newcommand{\taubin}[2]{\begin{bmatrix}{#1}\\ {#2}\end{bmatrix}_\tau}
\title{Blocks in ASEP with step-Bernoulli initial condition}

\begin{document}
	\maketitle
	
	\begin{abstract}
		This paper extends work by Tracy and Widom on blocks in the asymmetric simple exclusion process (ASEP) to the case of step-Bernoulli initial condition. We consider the probability that a particle at site $x$ is the beginning of a block of $L$ consecutive particles at time $t$ in ASEP with step-Bernoulli initial condition. A Fredholm determinant representation for this probability is derived, and the asymptotics are computed for the KPZ regime.
	\end{abstract}

	\section{Introduction}
	The \textit{asymmetric simple exclusion process} (ASEP) is a stochastic process on the integer lattice $\Z$ where each particle waits an exponential time, then jumps to the left with probability $p$ and to the right with probability $q=1-p$, unless the site to which it would jump is occupied in which case the particle remains where it is. \\
	
	In the case of step initial condition, where at time zero positive integer sites are occupied and all other sites are unoccupied, a formula for the distribution of the $m$th particle from the left \cite{asymptotics} was the starting point for the one-point probability distribution of the height function for the Kardar-Parisi-Zhang (KPZ) equation with narrow wedge initial conditions \cite{acq,ss}.\\
	
	In \cite{blocks}, a formula for the probability that a block of $L$ particles starts at site $x$ at time $t$ was derived; and in \cite{blockasymptotics}, the asymptotics (as $x,m\to\infty$ in the KPZ regime) of this probability were computed for step initial condition. In \cite{stepbernoulli}, a formula for the probability that the $m$th particle is in site $x$ at time $t$ was derived for step-Bernoulli initial condition, where positive integer sites are occupied with probability $\rho$ and the other sites are unoccupied. Asymptotics of this system (as $x,m,t\to\infty$) were computed in the KPZ regime. Here, we combine the two cases above and consider the probability that the $m$th particle is at the beginning of a block of length $L$ at time $t$ starting from step-Bernoulli initial condition, and then compute the asymptotics. Formally, this is the probability $\bbP_{L,\rho}(x,m,t)$ of the event
	$$
	x_m(t)=x, x_{m+1}(t)=x+1,\dots,x_{m+L-1}(t)=x+L-1
	$$
	starting from step-Bernoulli initial condition with parameter 
	$\rho$, where $x_m(t)$ is the position of the $m$th particle at time $t$.
	\\

	We work under the assumption that $q\neq0$ and we define $\tau=p/q$, $\g=q-p>0$. To state our first theorem, we must introduce some notation. Let $\z,\z'\in \bbC$ and define
	
	\begin{equation*}
	U(\z,\z')=\frac{p+q\z\z'-\z}{\z'-\z}
	\end{equation*}
	$$
	\epsilon(\z)=p\z^{-1}+q\z-1.
	$$
	
	Let $\cC_R$ be the circle centered at the origin with radius $R$ oriented counterclockwise and let $z\in\bbC^L$.
	Let $K_x$, $K_{x,\rho}$ and $K_{L,x,\rho}(z)$\footnote{We use the same notation for the operator as its kernel.} be the integral operators acting on $\cC_R$\footnote{When we say an operator acts on a subset of $\bbC$, we mean that it acts on $L^2$ functions on that set.} with kernels
	$$
	K_x(\x,\x')=\frac{\x^x e^{\e(\x)t}}{p+q\x\x'-\x},
	$$
	$$
	K_{x,\rho}(\xi,\xi')
	=
	qK_x(\x,\x') \frac{\rho(\x-\tau)}{\x-1+\rho(1-\tau)},
	$$
	and
	$$
	K_{L,x,\rho}(\x,\x';z)=K_{x,\rho}(\x,\x') \prod_{j=1}^{L}U(z_j,\x)
	$$
	respectively. Here $R$ is chosen to be large enough so that $\cC_R$ contains all the singularities of the integrand. 
	
	Let $(\la;\tau)_m=\prod_{j=0}^{m-1}(1-\la\tau^{j})$ be the $\tau$-Pochhammer symbol.
	
	All contour integrals are to be given a factor of $1/2\pi i$. The empty product is taken to be 1. 
	\\
	
	\noindent\textbf{Theorem 1}  \textit{Assume} $p,q\neq 0.$ \textit{Then}
	\begin{flalign*}
	\bbP_{L,\rho}(x,m,t)=\notag
	&(-1)^{L-1}p^{L(L+1)/2}\tau^{(1-m)(L-1)}\notag\\
	&\times \int_{\Gamma_{0,\tau}} \dots   \int_{\Gamma_{0,\tau}} \frac{1}{z_1^L(qz_1-p)z_2^{L-1}(qz_2-p)\cdots z_L(qz_L-p)}\notag\\
	&\times \left[ \int \frac{ \det(I-p^{-L}\lambda K_{L,x+L-1,\rho}(z))} {(\lambda;\tau)_m}\frac{d\lambda}{\lambda^L}\right] \prod_{1\leq i<j\leq L} \frac{1}{U(z_j,z_i)} dz_L\cdots dz_1.
	\end{flalign*}
	\textit{
	The integration with respect to $\lambda$ is over a contour containing the singularities of the integrand at $\la=\tau^{-j}$ for $j=0,\dots,m-1$, and at $\la=0$. The contours $\Gamma_{0,\tau}$ are small circles about the points $z=0$ and $z=\tau$ such that the $z_i$ contour is well inside the $z_{i-1}$ contour.}\\

	Notice when $\rho=1$ this is Theorem 3 of \cite{blocks}. Also, notice that when $L=1$ the $z_1$ integral is equal to the sum of the residues at 0 and $\tau$, which gives us
	$$
	\bbP_{1,\rho}(x,m,t)
	=
	\int \frac{\det(I-\la K_{x,\rho})-\det(I-\la K_{x-1,\rho})}{(\la;\tau)_m}
	\frac{d\la}{\la}.
	$$
	This implies
	$$
	\bbP_\rho(x_m(t)\leq x)
	=
	\int \frac{\det(I-\la K_{x,\rho})}{(\la;\tau)_m}
	\frac{d\la}{\la},
	$$
	which is Theorem 1 of \cite{stepbernoulli}.
	\\

	We assume throughout that $\s>0$ and we define
	$$
	m=\s t,
	\hspace{7mm}
	c_1=-1+2\sqrt{\s},
	\hspace{7mm}
	c_2=\s^{-1/6}(1-\sqrt{\s})^{2/3},
	$$
	$$
	c_1'=\rho^{-1}\s+\rho-1,
	\hspace{7mm}
	c_2'=\rho^{-1}((1-\rho)(\s-\rho^2))^{1/2}.
	$$
	
	Let $G(s)$ be the Gaussian distribution defined by
	$$
	G(s)=\frac{1}{\sqrt{2\pi}}\int_{-\infty}^s e^{-z^2/2}dz.
	$$
	Let
	$$
	K_{\Airy}(x,y)=\int_0^\infty \Ai(x+z)\Ai(y+z)dz,
	$$
	$$
	F_2(s)=\det(I-K_{\Airy}) \-\ \text{ on } (s,\infty),
	$$ 
	and let $F_1(s)^2$ be the given by the analogous determinant where $K_{\Airy}(x,y)$ is replaced by
	$$
	K_{\Airy}(x,y)+\Ai(x)\int_{-\infty}^y \Ai(z)dz.
	$$
	Although $F_1(s)^2$ is the square of a distribution from random matrix theory (see \cite{fs}), it arises here instead as the Fredholm determinant of the rank one perturbation of the Airy kernel.\\
	
	\noindent\textbf{Theorem 2}\\
	\textit{In the limit as} $t\to\infty,$
	\[
		\bbP_{L,\rho}(x,m,t/\g)=
		\begin{cases}
		c_2^{-1}\s^{(L-1)/2}F_2'(s)t^{-1/3}+o(t^{-1/3})
		&\text{ when } \s<\rho^2,
		\\
		c_2^{-1}\sigma^{(L-1)/2} \frac{d}{ds}\left(F_1(s)^2\right)t^{-1/3}+o(t^{-1/3})
		&\text{ when } \s=\rho^2,\-\ \rho<1,
		\\
		c'^{-1}_2\s^{(L-1)/2}G'(s)t^{-1/2}+o(t^{-1/2})
		&\text{ when } \s>\rho^2,\-\ \rho<1.
		\end{cases}
	\]
	\textit{Where in the first two cases
	$s=(x-c_1t)c_2^{-1}t^{-1/3}$ and in the third case $s=(x-c_1't)c_2'^{-1}t^{-1/2}$. These hold uniformly for $s$ in a bounded set, and uniformly for $\s$ in a compact subset of its domain.}\\

 	Notice that when $\s<\rho^2$ and $\rho=1$ we recover the result of \cite{blockasymptotics}. \\
	
	\noindent\textbf{Corollary} \\
	\textit{For all $0<\rho\leq1$, given that the $m$th particle from the left is in position $x$ at time $t$ the conditional probability that it is the beginning of a block of length $L$ converges to $\s^{(L-1)/2}$ as $t\to\infty$.}\\
	
	This follows from the fact that the conditional probability is equal to.
	$$
	\frac{\bbP_{L,\rho}(x,m,t/\g)}{\bbP_{1,\rho}(x,m,t/\g)}=\s^{(L-1)/2}+o(1).
	$$
	\\
	
	In Section \ref{freddet}, we prove Theorem 1. In Section \ref{Jker} we derive an alternate formula for $\bbP_{L,\rho}(x,m,t)$ that will be used in the proof of Theorem 2. In Section \ref{asy} we prove Theorem 2 in three cases. The first two cases of Theorem 2 are highly similar to \cite{blockasymptotics} and so the details are left to Appendices \ref{case1} and \ref{case2}. Appendix \ref{appwint} contains a result of Harold Widom \cite{widom} which is used in Section \ref{asy3}.
	\section{Proof of Theorem 1} \label{freddet}
	
	We begin by introducing notation. 
	Let
	$$
	f_L(\xi)=p^{L(L+1)/2-LN}\prod_i\xi_i^L
	\int_{\G_\xi}\cdots\int_{\G_\xi}
	\phi_L(z_1,\dots,z_L;\xi)
	dz_1\cdots dz_L,
	$$
	and
	$$
	I_L(x,k,\x)=\prod_{1\leq i<j\leq L}\frac{1}{U(\x_i,\x_j)} \prod_{i=1}^L\frac{1}{1-\x_i}f_L(\x) \prod_{i=1}^L\left(\x_i^{x-1}e^{\e(\x_i)t}\right).
	$$
	Here the $\G_\xi$ are simple closed curves containing the points $\xi_j$ but no other singularities of the integrand $\ph_L$, which is defined by
	$$
	\phi_L(z_1,\dots,z_L;\xi)
	=
	\frac{
			\prod_{1\leq j\leq N}
			U(z_1,\xi_j)U(z_2,\xi_j)\cdots U(z_L,\xi_j)
		}
		{
			z_1^L(qz_1-p)z_2^{L-1}(qz_2-p)
			\cdots
			z_L(qz_L-p)
		}
	\prod_{1\leq i<j\leq L} \frac{1}{U(z_j,z_i)}.
	$$
	as in \cite{blocks}. Let
	$$
	I_L(x,S,\x)=I_L(x,k,\x)\prod\xi_i^{-s_i}
	$$
	where $S=\{s_1,\dots s_k \}\subset \{1,\dots,N\}.$ For two sets $U,V\subset \Z$, define
	$$
	\s(U,V)=\#\{(u,v):u\in U, v\in V\}.
	$$
	Let $\bbP_{L,Y}(x,m,t)$ by the probability that the $mth$ particle from the left is at site $x$ at time $t$ starting from the initial conditions where the sites $Y$ are occupied and all other sites are unoccupied.
	Theorem 2 of \cite{blocks}, after a little algebra, is
	\begin{equation}
	\bbP_{L,Y}(x,m,t)=\sum_{k\geq m+L-1}c_{m,k,L}
	\sum_{\substack{S\subset Y \\ |S|=k}}\tau^{\sigma(S,Y)}
	\int_{|\xi_i|=R}I_L(x,S,\xi)d^{|S|}\xi,
	\label{PY}
	\end{equation}
	where $R$ is large enough for the poles of $I_L$ to be contained inside the contour,
	$$
	c_{m,k,L}=q^{k(k-1)/2}(-1)^{m+1}\tau^{m(m-1)/2} \tau^{-km}\taubin{k-L}{m-1},	
	$$
	and
	$$
	\taubin{n}{m}=\prod_{i=1}^{m}\frac{1-\tau^{n-i-1}}{1-\tau^{i-1}}
	$$
	is the $\tau$-binomial coefficient. \\
	
	To find the probability for step-Bernoulli initial condition, we first take a weighted average over all initial configurations $Y\subset\{1,...N\}$. As in \cite{stepbernoulli}, the only factor in (\ref{PY}) that depends on $Y$ is $\tau^{\sigma(S,Y)}$ and the probability of the initial condition $Y$ is given by $\rho^{|Y|}(1-\rho)^{N-|Y|}$. From \cite{stepbernoulli} we have that
	\begin{equation}
	\sum_{S\subset Y \subset [1,N]}\rho^{|Y|}(1-\rho)^{N-|Y|}\tau^{\sigma(S,Y)}
	=
	\tau^{k(k+1)/2}\rho^k\prod_{i}(1-\rho+\tau^{k-i+1}\rho)^{t_i}
	\label{sumoverY}
	\end{equation}
	where we define $t_i=s_i-s_{i-1}-1$. 
	
	Continuing to follow \cite{stepbernoulli}, we notice that the only factor in (\ref{PY}) which depends on the elements of $S$ (rather than just the size of $S$) is
	$$
	\prod\xi_i^{-s_i}=\prod\xi_i^{-(t_1+\cdots+t_i+i)}= \prod\xi_i^{-i} \cdot\prod(\xi_i\xi_{i+1}\cdots\xi_k)^{-t_i}.
	$$
	Multiplying this by (\ref{sumoverY}) and summing over all $S\subset\{1,...N\}$ we obtain
	\begin{equation*}
		\tau^{k(k+1)/2}\rho^k\prod\x_i^{-i}\sum_{\sum t_i\leq N-k} \prod_i\left(\frac{1-\rho+\tau^{k-i+1}\rho}{\x_i\cdots\x_k}\right)^{t_i}.
	\end{equation*}
	Letting $N\to\infty$ then gives us
	\begin{equation*}
		\tau^{k(k+1)/2}\rho^k\prod_i \frac{1}{\x_i\cdots\x_k-1+\rho-\tau^{k-i+1}\rho}.
	\end{equation*}
	And so we may conclude that
	\begin{flalign*}
		\bbP_{L,\rho}(x,m,t) = &\sum_{k\geq m+L-1}\tau^{k(k+1)/2} \rho^k c_{m,k,L}\\
		&\times
		\int_{|\xi_i|=R} I_L(x,k,\xi)\prod_i \frac{1}{\xi_i\cdots\xi_k-1+\rho-\tau^{k-i+1}\rho}\prod_i d\xi_i
	\end{flalign*}
	
	Notice that $f_L$ is symmetric in the $\x_i$ and so applying combinatorial identity (9) from \cite{stepbernoulli} gives us
	\begin{flalign*}
		\bbP_{L,\rho}(x,m,t) = 
		&\sum_{k\geq m+L-1}\frac{1}{k!} q^{k(k-1)/2}\tau^{k(k+1)/2}c_{m,k,L}\\
		&\times\int_{|\xi_i|=R}\prod_{i\neq j} \frac{1}{U(\x_i,\x_j)}f_L(\x)\\
		&\times\prod_i\frac{\rho}{\x_i-1+\rho(1-\tau)} \prod_i\frac{\x_i^{x-1}e^{\e(\x_i)t}}{1-\xi_i}d\x_i.
	\end{flalign*}

	Applying the alternative expression of $f_L$ from \cite{blocks}, the above integrand becomes
	$$
	(-1)^Lp^{L(L-1)/2-kL} \int_{\Gamma_{0,\tau}} \dots \int_{\Gamma_{0,\tau}} \frac{1}{z_1^L(qz_1-p)z_2^{L-1}(qz_2-p)\cdots z_L(qz_L-p)} \prod_{i<j} \frac{1}{U(z_j,z_i)}	
	$$
	\begin{equation}{\label{bracket}}
	\times \left[\prod_{i\neq j}\frac{1}{U(\x_i,\x_j)} \prod_i\frac{U(z_1,\x_i)U(z_2,\x_i)\cdots U(z_L,\x_i)}{(1-\x_i)(q\x_i-p)} \prod_i \frac{q\rho(\x_i-\tau)}{\x_i+\rho(1-\tau)}\prod_i \x_i^{x+L-1}e^{\e(\x_i)t} \right] dz_L\cdots dz_1 
	\end{equation}
	Where the contour $\Gamma_{0,\tau}$ consists of tiny circles around the points $0,\tau$ with the $z_i$ contours well inside the $z_{i-1}$ contours.

	From equation (3) of \cite{fredholm} we have 
	$$
	\det(K_x(\x_i,\x_j))_{i,j\leq k}=(-1)^k(pq)^{k(k-1)/2} \prod_{i\neq j}\frac{1}{U(\x_i,\x_j)}\prod_i\frac{1}{(1-\x_i)(q\x_i-p)} \prod_i \x_i^x e^{\e(\x_i)t}
	$$
	
	Thus, the expression in brackets in (\ref{bracket}) is equal to
	
	\begin{equation*}
		(-1)^k(pq)^{-k(k-1)/2}\det(K_{L,x+L-1,\rho}(\x_i,\x_j;z))_{i,j\leq k}
	\end{equation*}
	
	Notice that
	\begin{equation*}
		\frac{(-1)^k}{k!} \int_{|\x_i|=R} \cdots \int_{|\x_i|=R} \det(K_{L,x+L-1,\rho}(\x_i,\x_j;z))_{i,j\leq k} d\x_i\cdots d\x_k
	\end{equation*}
	\begin{equation*}
		=\int\frac{\det(I-\lambda K_{L,x+L-1,\rho}(z))} {\lambda^{k+1}} d\lambda 
		=p^{kL}\int\frac{\det(I-\lambda p^{-L}K_{L,x+L-1,\rho}(z))} {\lambda^{k+1}} d\lambda
	\end{equation*}
	the $k$th coefficient in the Fredholm expansion of $\det(I-\lambda K_{L,x+L-1,\rho}(z))$.

	Next we sum over $k$ and bring the sum inside the integrals. Considering terms involving $k$, we have
	\begin{equation*}
		\sum_{k\geq L+m-1} q^{k(k-1)/2} \tau^{k(k+1)/2}  (pq)^{-k(k-1)/2} \lambda^{-k} q^{k(k-1)/2} (-1)^{m+1}\tau^{m(m-1)/2}\tau^{-km} \taubin{k-L}{m-1}
	\end{equation*}
	\begin{equation}  
		=(-1)^{m+1}\lambda^{-L-m+1}\tau^{(1-m)(L-1)}\tau^{-m(m-1)/2} \sum_{k\geq0} \tau^{(1-m)k}\lambda^{-k} \taubin{k+m-1}{k}
		\label{kterms}
	\end{equation}
	
	By the $\tau$-binomial theorem, for $|\tau^{1-m}|<|\lambda|$ we have that (\ref{kterms}) is equal to
	\begin{equation*}
		(-1)^m\tau^{-m(m-1)/2}\la^{-m}\sum_{k\geq0}\left(\frac{\tau^{1-m}}{\la}\right)^k\taubin{k+m-1}{k}
	\end{equation*}
	\begin{equation*}
		= (-1)^m\tau^{-m(m-1)/2}\la^{-m}\prod_{j=0}^{m-1} \left(1-\frac{\tau^{1-m}}{\la}\tau^j\right)^{-1} 
		= \frac{1}{(\lambda;\tau)_m}
	\end{equation*}
	
	Thus Theorem 1 holds.

	\section{Replacing the operators $K_{L,x,\rho}(z)$ by operators $J_{L,x,m,\rho}(w)$} \label{Jker}
	First we must introduce notation. Let
	$$
	\ph_\infty(\eta)= (1-\eta)^{-x-L+1}e^{\frac{\eta}{1-\eta}t} \prod_{n=0}^\infty\frac{1}{1+\tau^n\alpha\eta},
	$$
	where, as in \cite{stepbernoulli}, $\alpha=\frac{1-\rho}{\rho}$. Let
	
	$$
	f(\mu,z) =\sum_{k\in\Z}\frac{\tau^k}{1-\tau^k\mu}z^k, 
	\-\ 
	\text{ and }
	\-\
	V(\z,\eta';w)=\frac{w\z-\tau}{w\eta'-\tau}.
	$$
	Let $J_{L,x,m,\rho}(w)$ be the operator with kernel
	\begin{equation}
	J_{L,x,m,\rho}(\eta,\eta';w)= \int\frac{\ph_\infty(\z)}{\ph_\infty(\eta)} \frac{\z^{m-L}}{(\eta')^{m-L+1}} \frac{f(\mu,\z/\eta')}{\z-\eta} \prod_{j=1}^L V(\z,\eta';w_j)d\z.
	\label{J}
	\end{equation}
	\\
	When $\rho>1/2$ the singularities $-(\tau^n\alpha)^{-1}$ lie outside the unit circle and we take the $\z$-contour to be a circle slightly larger than the unit circle. In this case the operator $J_{L,x,m,\rho}$ acts on a circle slightly smaller than the unit circle. When $\rho\leq1/2$ the $\z$-contour is the circle with diameter $[-\alpha^{-1}+\delta,1+\delta]$ and $J_{L,x,m,\rho}$ acts on the circle with diameter $[-\alpha^{-1}+2\delta,1-\delta]$.\\
	
	\noindent\textbf{Lemma} \textit{Let the contours $\Gamma_{0,\tau}$ be defined as in Theorem 1.  Then}
	\begin{equation*} 
	\bbP_{L,\rho}(x,m,t)=-\tau^{-(L^2-5L+2)/2}\int_{\Gamma_{0,\tau}}\cdots\int_{\Gamma_{0,\tau}}\prod_{j=1}^{L}\frac{(w_j-1)^{L-j}}{w_j(w_j-\tau)^{L-j+1}}\prod_{i<j}\frac{w_j- w_i}{w_j- \tau w_i}
	\end{equation*}
	\begin{equation}
	\times\int\left[ (\tau^L\mu;\tau)_\infty \det(I+\mu J_{L,x,m,\rho}(w)) \frac{d\mu}{\mu^L}\right]dw_L\cdots dw_1,
	\label{J-rep}
	\end{equation}
	\textit{where the order of integration of the $w_i$ is as indicated.}\\
	
	The only difference between this lemma and Section 2 of \cite{blockasymptotics} is the infinite product in $\ph_\infty$. When $\rho>1/2$, the singularities in this product are all greater than $\tau^{-1}$	so the argument from \cite{blockasymptotics} goes over unchanged. 
	
	The right side of Theorem 1 is an analytic function of $\rho$ on $(0,1)$. The right side of (\ref{J-rep}) can be extended to an analytic function of $\rho$ on $(0,1)$ by varying the $\z$, $\eta$, and $w_i$-contours such that the singularities $-(\tau^n\alpha)^{-1}$ lie outside the $\z$-contour, the $\eta$-contour remains inside the $\z$-contour, and the $w_i$-contour remains inside the $\eta$ contour. The expressions agree for $\rho\in(1/2,1]$ and so agree for all $\rho$. 
	
	\section{Proof of Theorem 2} \label{asy}
	Theorem 2 naturally breaks into three cases. In all cases, we assume $m=\sigma t$. In the first two cases, we assume $x=c_1t+sc_2t^{1/3}$
	and, as in \cite{blockasymptotics}, compute the asymptotics using the 
	substitutions
	\begin{equation}
		\eta\to \xi+c_2^{-1}t^{-1/3}\eta, \hspace{5mm}
		\eta'\to \xi+c_2^{-1}t^{-1/3}\eta', \hspace{5mm} 
		\zeta\to \xi+c_2^{-1}t^{-1/3}\zeta.	
		\label{sub}
	\end{equation}

	\subsection{Case 1: $\s<\rho^2$}
	When $\sigma<\rho^2$ the argument from \cite{blockasymptotics} goes through unchanged since the pole structure of the integrand is the same and the value of $\ph_\infty(\zeta)/\ph_\infty(\eta')$ tends to 1 near the saddle point. Thus, we may conclude that, as $t\to\infty$,
	\begin{equation*}
	\bbP_{L,\rho}(x,m,t/\gamma)=c_2^{-1}\sigma^{(L-1)/2} F_2'(s)t^{-1/3}+o(t^{-1/3}).
	\end{equation*}
	This is the first part of Theorem 2. For a more detailed argument see Appendix \ref{case1}.

	\subsection{Case 2: $\s=\rho^2$}
	In this case, as in \cite{blockasymptotics}, we still have that the operator $\mu J_{x,L,m,\rho}$ has the same Fredholm determinant as the sum of
	\begin{equation*}
	J^{(0)}+o(1)
	\end{equation*}
	and
	$$
	J^{(1)}\sum_{j=1}^L\frac{w_j}{w_j\xi-\tau}c_3^{-1}t^{-1/3}+o(t^{-1/3})
	$$
	where $o(1)$ and $o(t^{-1/3})$ denote operators whose trace norms are $o(1)$ and $o(t^{-1/3})$ respectively. However, $J^{(0)}$ and $J^{(1)}$ have kernels
	\begin{equation*}
		J^{(0)}(\eta,\eta')=\int_{\Gamma_\z} \frac{e^{\z^3/3+s\z+(\eta')^3/3-s\eta'}}{(\z-\eta)(\eta'-\z)} \frac{\eta}{\z}d\z
	\end{equation*}
	and
	\begin{equation*}
	J^{(1)}(\eta,\eta')=-\int_{\Gamma_\z} \frac{e^{\z^3/3+s\z+(\eta')^3/3-s\eta'}}{\z-\eta} \frac{\eta}{\z}d\z.
	\end{equation*}
	This is because, after the substitutions (\ref{sub}), 
	\begin{equation}
	\prod_{j=1}^LV(\z,\eta';w_j)\to 1+(\z-\eta)\sum_{j=1}^L\frac{w_j}{w_j\xi-\tau}c_3^{-1}t^{-1/3}+E(\z,\eta';w),
	\label{sum}
	\end{equation}
	where $E(\z,\eta';w)$ is a polynomial in $\z-\eta'$ with $O(t^{-2/3})$ coefficients.
	If we distribute the other factors in our kernel through this sum, we can consider $J_{L,x,m,\rho}$ as the sum of three operators. The operator corresponding to the additive 1 in (\ref{sum}) is $J^{(0)}+o(1)$, by the argument in \cite{stepbernoulli}. Although our kernel has $m-L$ where the kernel in \cite{stepbernoulli} has $m$, this results only in an $O(1/t)$ change in $\sigma$. The argument in \cite{stepbernoulli} begins by setting $\s=\rho^2+o(t^{-1/3})$ and so the $O(1/t)$ change makes no difference.
	Similarly, the operator corresponding to the second term in (\ref{sum}) is 
	$$
	J^{(1)}\sum_{j=1}^L\frac{w_j}{w_j\xi-\tau}c_3^{-1}t^{-1/3}+o(t^{-1/3}).
	$$
	 
	The operator corresponding to the error term $E(\z,\eta';w)$ has trace norm $O(t^{-2/3})$. 
	
	Since $J^{(1)}=\frac{d}{ds}J^{(0)}$ and $J^{(0)}$ is independent of $w$, the argument from \cite{blockasymptotics} can be followed word for word to conclude
	\begin{equation*}
		\bbP_{L,\rho}(x,m,t/\gamma)=c_2^{-1}\sigma^{(L-1)/2} \frac{d}{ds}\left(F_1(s)^2\right)t^{-1/3}+o(t^{-1/3}).
	\end{equation*}
	\\
	
	\subsection{Case 3: $\s>\rho^2$} \label{asy3}
	Define $\psi$ such that
	\begin{equation}
	e^{\psi(\z)}=(1-\z)^{-x-L+1}e^{\frac{\z}{1-\z}t}\z^{m-L}
	\label{psi}.
	\end{equation}
	
	In \cite{stepbernoulli} the determinant of the $J$-kernel was computed by deforming the $\eta$-contour to the component of the level curve $|e^{\psi(\eta)}/e^{\psi(-\alpha^{-1})}|=1-\delta$ with the saddle point outside and a small indentation to the left of $\eta=1$, and the $\z$-contour to the component of the level curve $|e^{\psi(\z)}/e^{\psi(-\alpha^{-1})}|=1-2\delta$ with the saddle point inside and a small indentation to the right of $\z=1$. This resulted in the norm of the operator represented by the new $\z$ integral being exponentially small.\\
	The extra factors found when considering blocks do not depend on $t$, and so they only change the norm of the operator by $O(1)$. Thus, if we change the contours of $J$ as above, we pick up the residue
	$$
	\frac{e^{\psi(-\alpha^{-1})}}{e^{\psi(\eta)}}\eta^{-1} f(\mu,-(\alpha\eta)^{-1})\left(\frac{1-\eta}{1+\alpha^{-1}}\right)^{c_2'st^{1/2}} \prod_{n=1}^\infty \frac{1+\tau^n\alpha\eta}{1-\tau^n}\prod_{j=1}^{L}\frac{-\alpha^{-1}w-\tau}{w\eta-\tau}
	$$
	from the pole at $\z=-\alpha^{-1}$. Therefore, with exponentially small error in $t$, $\det(I+\mu J)$ equals

	\begin{equation*}
	1+\mu\int
	\frac{e^{\psi(-\alpha^{-1})}}{e^{\psi(\eta)}}\eta^{-1} f(\mu,-(\alpha\eta)^{-1})\left(\frac{1-\eta}{1+\alpha^{-1}}\right)^{c_2'st^{1/2}} \prod_{n=1}^\infty \frac{1+\tau^n\alpha\eta}{1-\tau^n}\prod_{j=1}^{L}\frac{-\alpha^{-1}w-\tau}{w\eta-\tau}
	d\eta.
	\end{equation*}
	We can ignore the first summand (which is 1) because, by appendix $B$ of \cite{blocks}, 
	$$
	\int_{\Gamma_{0,\tau}}\cdots\int_{\Gamma_{0,\tau}}
	\prod_{j=1}^{L}\frac{(w_j-1)^{L-j}}{w_j(w_j-\tau)^{L-j+1}}\prod_{i<j}\frac{w_j-w_i}{w_j-\tau w_j}
	dw_L\cdots dw_1=0.
	$$
	
	For the second summand, exchanging the order of integration and applying the result of appendix \ref{appwint} gives us
	
	$$
	\int_{\Gamma_{0,\tau}}\cdots\int_{\Gamma_{0,\tau}}\prod_{j=1}^{L}\frac{(w_j-1)^{L-j}}{w_j(w_j-\tau)^{L-j+1}}\prod_{i<j}\frac{w_j-w_i}{w_j-\tau w_j}
	\mu\int\frac{e^{\psi(-\alpha^{-1})}}{e^{\psi(\eta)}}\eta^{-1} f(\mu,-(\alpha\eta)^{-1})
	$$
	$$
	\times	\left(\frac{1-\eta}{1+\alpha^{-1}}\right)^{c_2'st^{1/2}} \prod_{n=1}^\infty \frac{1+\tau^n\alpha\eta}{1-\tau^n}\prod_{j=1}^{L}\frac{-\alpha^{-1}w-\tau}{w\eta-\tau}d\eta
	dw_L\cdots dw_1
	$$
	$$
	=(-1)^{L-1}\mu\int\frac{\psi(-\alpha^{-1})}{\psi(\eta)}\eta^{-1} f\left(\mu,\frac{-1}{\alpha\eta}\right)
	\left(\frac{1-\eta}{1+\alpha^{-1}}\right)^{c_2'st^{1/2}}
	$$
	\begin{equation}
	\times\prod_{n=1}^\infty\frac{1+\tau^n\alpha\eta}{1-\tau^n}
	\frac{\prod_{J=0}^{L-1}(\alpha\eta+\tau^j)}{\alpha^L(\eta-1)^L\tau^{L^2}}
	d\eta,
	\label{eta}
	\end{equation}
	where
	$x=c'_1t+c'_2st^{1/2}$.
	We continue to follow the argument from \cite{stepbernoulli}, evaluating the integral by using the steepest descent curve with saddle point $-\alpha^{-1}$. The extra factors give us
	$$
	\eta^{-1} f\left(\mu,\frac{-1}{\alpha\eta}\right)\frac{\prod_{J=0}^{L-1}(\alpha\eta+\tau^j)}{\alpha^L(\eta-1)^L\tau^{L^2}}
	=
	\frac{(\tau-1)(\tau^2-1)\cdots(\tau^{L-1}-1)}{\mu\alpha^{L-1}(-\alpha^{-1}-1)^L\tau^{L^2}}+O\left(\eta+\alpha^{-1}\right).
	$$
	Near the saddle point, we have
	$$
	\left(\frac{1-\eta}{1+\alpha^{-1}}\right)^{c_2'st^{1/2}}
	=
	e^{-[(1-\rho)c_2's(\eta+\alpha^{-1})+O((\eta+\alpha^{-1})^2)]t^{1/2}}
	$$
	$$
	=\frac{d}{ds}\frac{t^{-1/2}}{(\rho-1)c_2'(\eta+\alpha^{-1})}e^{-[(1-\rho)c_2's(\eta+\alpha^{-1})+O((\eta+\alpha^{-1})^2)]t^{1/2}},
	$$
	the infinite product in the integrand is $1+O(\eta+\alpha^{-1})$, and 
	$$
	\frac{e^{\psi(-\alpha^{-1})}}{e^{\psi(\eta)}}
	=
	e^{[\rho^{-1}(1-\rho)^3(\s-\rho^2)(\eta+\alpha^{-1})^2/2+O((\eta+\alpha^{-1})^3)]t}
	$$
	Thus, for large $t$, the integral (\ref{eta}) is
	$$
	\frac{(\tau-1)\cdots(\tau^{L-1}-1)}{(-\alpha)^{L-1}(-\alpha^{-1}-1)^L\tau^{L^2}}
	\frac{t^{-1/2}}{(\rho-1)c_2'}
	\frac{d}{ds}
	\int_{-i\infty+0}^{i\infty+0}
	e^{\rho^{-1}(1-\rho)^3(\s-\rho^2)t\eta^2/2-(1-\rho)c_2'st^{1/2}\eta}
	\frac{d\eta}{\eta+\alpha^{-1}}+o(t^{-1/2})
	$$

	$$
	=\frac{(\tau-1)\cdots(\tau^{L-1}-1)}{(-\alpha)^{L-1}(-\alpha^{-1}-1)^L\tau^{L^2}}
	\frac{t^{-1/2}}{(\rho-1)c_2'}
	\frac{d}{ds}
	\int_{-i\infty+0}^{i\infty+0}
	e^{\eta^2/2-s\eta}
	\frac{d\eta}{\eta}+o(t^{-1/2})
	$$
	\begin{equation}
	=-\frac{(\tau-1)(\tau^2-1)\cdots(\tau^{L-1}-1)}{(-\alpha)^{L-1}(-\alpha^{-1}-1)^L\tau^{L^2}}
	\frac{t^{-1/2}}{(1-\rho)c'_2}G'(s)+o(t^{-1/2}).
	\label{G}
	\end{equation}
	Notice the factor $(-\alpha)^{-(L-1)}(-\alpha^{-1}-1)^{-L}=\sqrt{\s}^{L-1}(\xi-1)^{-1}$, from the definition of $\s$ and the fact that $\xi=-\alpha^{-1}$.
	From section III of \cite{blocks} we have
	$$
	\int (\tau^L\mu,\tau)_\infty\frac{d\mu}{\mu^L}
	=\frac{\tau^{(L-1)(3L-2)/2}}{(\tau-1)\cdots(\tau^{L-1}-1)}.
	$$
	Furthermore $\xi-1=1/(\rho-1)$ so $(\xi-1)(1-\rho)=-1$.
	Combining this with (\ref{G}) and (\ref{J-rep}) gives us that 
	$$
	\bbP_{L,\rho}(x,m,t/\g)
	=
	t^{-1/2}c'^{-1}_2\sqrt{\s}^{L-1}G'(s)+o(t^{-1/2})
	$$
	as $t\to\infty.$
	\\
	
	\begin{appendices}
		\section{$w$-integrals}
		\label{appwint}
		
		The following proof is due to Harold Widom \cite{widom}.
		\\
		
		All integrals in this appendix are over $\Gamma_{0,\tau}$. Let
		\begin{equation}
		\Phi_L(w_1,\dots,w_L)
		=
		\prod_{j\geq1}\frac{(w_j-1)^{L-j}}{w_j(w_j-\tau)^{L-j+1}}
		\prod_{i<j}\frac{w_i-w_j}{\tau w_i-w_j}
		\prod_{j\geq1}\frac{\ph_j(w_j)}{\theta-w_j}
		\label{Phi}
		\end{equation}
		where $\ph_j$ are affine transformations. We will show that 
		\begin{equation}
		\int\cdots\int \Phi_L(w_1,\dots,w_L) dw_L\cdots dw_1
		=
		\frac{1}{\theta^L(\theta-\tau)^L\tau^{L^2-L}}
		\prod_{j\geq1}\ph_j(\tau^{j-1}\theta)
		\label{wint}
		\end{equation}
		by considering the more general integral
		\begin{equation}
		V_{L,k}
		:=
		\int\cdots\int \prod_{j\geq1}\frac{\theta-w_j}{\tau^k\theta-w_j}\Phi_L(w_1,\dots,w_L) dw_L\cdots dw_1.
		\end{equation}
		The integrand may be rewritten as 
		\begin{equation}
		\frac{(w_1-1)^{L-1}}{w_1(w_1-\tau)^L}\frac{\ph_1(w_1)}{\tau^k\theta-w_1}G_L(w_1,\dots,w_L)\psi(w_2,\dots w_L)
		\label{wintegrand}
		\end{equation}
		where
		\begin{equation*}
		G_L(w_1,\dots w_L)
		=
		\prod_{j>1}\frac{(w_j-1)^{L-j}}{w_j(w_j-\tau)^{L-j+1}}
		\prod_{i<j}\frac{w_i-w_j}{\tau w_i-w_j}
		\end{equation*}
		and
		$$
		\psi(w_2,\dots,w_L)=\prod_{j\geq2}\frac{\ph_j(w_j)}{\tau^k \theta -w_j}.
		$$
		$\psi$ is analytic in a neighborhood of $\{0,\tau\}^{L-1}$ provided that $\theta\neq0, \tau^{1-k}$, so by appendix B of \cite{blocks}
		$$
		\int\cdots\int G_L(w_1,\dots,w_L)\psi(w_2,\dots w_L) dw_L\cdots dw_2
		$$
		is an analytic function of $w_1$ outside of $\{0,\tau\}$ except for a pole of order at most $L-1$ at 1. This pole is canceled by the first factor in (\ref{wintegrand}), so if we expand the contour we encounter only the pole at $w_1=\tau^k\theta$
		
		Furthermore, the assumption that the $\ph_j$ are affine ensures that there is no contribution at infinity when we expand the contour. Thus, we can conclude that
		\begin{equation*}
		V_{L,k}=\frac{(\tau^k\theta-1)^{L-1}}{\theta\tau^{k+L}(\tau^{k-1}\theta-1)^L}\ph_1(\tau^k\theta)V_{L-1,k+1}.
		\end{equation*}
		Applying the above results to
		\begin{equation*}
		\Phi_{L-n}(w_{n+1},w_{n+2},\dots,w_L)
		:=
		\prod_{j>n}\frac{(w_j-1)^{L-j}}{w_j(w_j-\tau)^{L-j+1}}
		\prod_{i<j}\frac{w_i-w_j}{\tau w_i-w_j}
		\prod_{j>n}\frac{\ph_j(w_j)}{\theta-w_j}
		\end{equation*}
		gives us
		\begin{equation*}
		V_{L-n,k}=\frac{(\tau^k\theta-1)^{L-n-1}}{\theta\tau^{k+L-n}(\tau^{k-1}\theta-1)^{L-n}}\ph_{n+1}(\tau^k\theta)V_{L-n-1,k+1},
		\end{equation*}
		from which it easily follows that
		\begin{equation*}
		V_{L,k}
		=
		\frac{1}{\theta^L\tau^{kL+L^2}(\tau^{k-1}\theta-1)^L}
		\ph_1(\tau^k\theta)\cdots\ph_L(\tau^{k+L-1}\theta).
		\end{equation*}
		Setting $k=0$ in the above equation gives us (\ref{wint}).
		
		Recall the notation 
		\begin{equation*}
		F_L(w_1,\dots w_L)
		=
		\prod_{j\geq1}\frac{(w_j-1)^{L-j}}{w_j(w_j-\tau)^{L-j+1}}
		\prod_{i<j}\frac{w_i-w_j}{\tau w_i-w_j}.
		\end{equation*}
		For the $w$-integral in \cite{blocks}, setting
		\begin{equation*}
		\theta=\tau/\xi,
		\hspace{3mm}
		\ph_1(w)=-w/\xi,
		\hspace{3mm}
		\ph_2(w)=\cdots=\ph_L(w)=\theta-w
		\end{equation*}
		in (\ref{Phi}) gives us, after some algebra,
		\begin{equation*}
		\int\cdots\int
		F_L(w_1,\dots,w_2)
		\frac{w_1}{w_1\xi-\tau}
		dw_L\dots dw_1
		=
		-\frac{(1-\tau)(1-\tau^2)\cdots(1-\tau^{L-1})}{\tau^{L^2}}
		\frac{\xi^{L-1}}{(1-\xi)^L}.
		\end{equation*}
		
		For the integral used in the $\sigma>\rho^2$ case
		setting
		\begin{equation*}
		\theta=\tau/\eta,
		\hspace{3mm}
		\ph_1(w)=\cdots=\ph_L(w)=-(w\alpha^{-1}+\tau)/\eta
		\end{equation*}
		in (\ref{Phi}) gives us, after some algebra,
		\begin{equation*}
		\int\cdots\int
		F_L(w_1,\dots,w_2)
		\prod_{j\geq1}\frac{w\alpha^{-1}+\tau}{w_j\eta-\tau}
		dw_L\dots dw_1
		=
		\frac{\prod_{j\geq1}(\alpha\eta+\tau^{j-1})}{\alpha^L(\eta-1)^L\tau^{L^2}}.
		\end{equation*}
		
		\section{Asymptotics when $\s<\rho^2$}
		\label{case1}
		Here we assume $m-L=\s t$ and $x+L-1=c_1t+c_2st^{1/3}$, but the following argument holds for $m=\s t$ and $x=c_1t+c_2st^{1/3}$ since changing $m$ and $x$ by $-L$ and $L-1$ is equivalent to changing $\s$ and $s$ by $Lt^{-1}$ and $(1-L)t^{-1/3}$ respectively.
		
		First, we consider the expression
		\begin{equation*}
			\log\big((1-\z)^{-x-L+1}e^{\frac{\z}{1-\z}t}\z^{m-L}\big)=(x+L-1)\log(1-\z)+\frac{\z}{1-\z}t+(m-L)\log(\z).
		\end{equation*}
		Differentiating this expression provides the saddle point equation
		
		\begin{equation*}
			\frac{x+L-1}{1-\z}+\frac{t}{(1-\z)^2}+\frac{m-L}{\z}=0.
		\end{equation*}
		 The two saddle points coincide when
		 \begin{equation}
		 	m-L=\frac{(x+L-1+t)^2}{4t}.
		 	\label{m-L}
		 \end{equation}
		 If $m-L=\s t$ and $x+L-1=c_1 t$ then (\ref{m-L}) gives us
		 
		 \begin{equation*}
		 	\s=\frac{(c_1+1)^2}{4}
		 \end{equation*}
		 so $c_1=-1+2\sqrt{\s}$ (we take the positive root because $c_1$ should increase with $\s$). Thus, the saddle point is at 
		 \begin{equation*}
		 	\xi=-\frac{\sqrt{\s}}{1-\sqrt{\s}}.
		 \end{equation*}
		 Let $c_2=\s^{-1/6}(1-\sqrt{\s})^{2/3}$ and set $x+L-1=c_1t+c_2st^{1/3}$. Let $\psi$ be as in (\ref{psi}).
		 Taking the Taylor expansion of $\psi$ around the point $\z=\xi$, we have
		 \begin{equation}
		 	\psi(\z)=-c_3^3t(\z-\xi)^3/3+c_3st^{1/3}(\z-\xi)+O(t(\z-\xi)^4)+O(t^{1/3}(\z-\xi)^2).
		 	\label{cubic}
		 \end{equation}
		 where $c_3=\s^{-1/6}(1-\sqrt{\s})^{5/3}$
		 Define $\psi_0(\z)=(x+L-1)\log(1-\z)+\frac{\z}{1-\z}t+(m-L)\log(\z)$ and $\psi_1(\z)=\psi_0(\z)-\psi_0(\xi)$.
		 Lemma 5 of \cite{asymptotics} gives us the existence of contours $\G_\z$ and $\G_\eta$ with the following properties: \\
		 (i) The part of $\G_\eta$ in a neighborhood $N_\eta$ of $\eta=\xi$ is a pair of rays $\xi$ in the directions $\pm\pi/3$ and the part of $\G_\z$ in a neighborhood $N_\z$ of $\z=\xi$ is a pair of rays from $\xi-t^{-1/3}$ in the directions $\pm2\pi/3$.\\
		 (ii) For some $\delta>0$ we have $\text{Re}(\psi_1(\z))<-\delta$ on $\G_\z\setminus N_\z$ and $\text{Re}(\psi_1(\eta))>\delta$ on $\G_\eta$.\\
		 (iii) The circular $\eta$- and $\z$--contours for $J_{L,x,m,\rho}$ can be simultaneously deformed to $\G_\eta$ and $G_\z$ respectively, so that during the deformation the integrand in (\ref{J}) remains analytic in all variables (this requires $\sigma<\rho^2$).\\
		 By condition (iii) of the above lemma and Proposition 1 of \cite{asymptotics}, $\det J$ remains the same if $J$ acts on $\G_\eta$ and the integral (\ref{J}) is over $\G_\z$.\\
		 Aside from the $e^{\psi(\z)}/e^{\psi(\eta)}$ factor, the $J$-kernel is uniformly $O(t^{1/3})$. This is because $\z-\eta=t^{-1/3}$, and none of the remaining factors grow with $t$. Together with (ii), this gives us that when we restrict to $\z\in\G_\z\setminus N_\z$ and $\eta\in\G_\eta\setminus N_\eta$, $J$ has exponentially small trace norm. For $a<1/3$, we may further restrict $\eta$ and $\zeta$ to rays of length $t^{-a}$, since the kernel has trace norm $O(e^{\delta t^{1-3a}})$ outside of a $t^{-a}$-neighborhood of $\xi$ by (\ref{cubic}). 
		 
		 On these rays, we make the substitutions (\ref{sub}).
		 Each $V(\z,\eta';w)$ becomes $1+(\z-\eta')\frac{w}{w\xi-\tau}c_3^{-1}t^{-1/3}[1+O(\min(1,t^{-1/3}|\eta'|))]$,
		 the product $\prod V(\z,\eta';w_j)$ becomes (\ref{sum}), and so $J$ can be considered as $T_1+T_2+T_3$, corresponding to the terms $1$, $(\z-\eta)\sum\frac{w_j}{w_j\xi-\tau}c_3^{-1}t^{-1/3}$, and $E(\z,\eta';w)$ in (\ref{sum}).
		 
		 $T_1$ can be written as the product $A_1B_1$ where, before substitution, $A_1:L^2(\G_\z)\to L^2(\G_\eta)$ and $B_1:L^2(\G_\eta)\to L^2(\G_\z)$ have
		 kernels
		 \begin{equation}
		 	A_1(\eta,\z)=\frac{e^{\psi(\z)}}{\z-\eta}\prod_{n=0}^\infty \frac{1+\tau^n\alpha\eta}{1+\tau^n\alpha\z}, \hspace{4 mm} 
		 	B_1(\z,\eta')=\frac{\mu f(\mu,\z/\eta)}{\eta e^{\psi(\eta)}}
		 	\label{AB}
		 \end{equation}
		 respectively. After substitution, the factor $1/(\z-\eta)$ is unchanged.
		 The product in the kernel of $A_1$ becomes
		 \begin{equation}
		 	\prod_{n=0}^\infty \frac{c_3t^{1/3}(\alpha^{-1}+\tau^n \xi)+\tau^n\eta}{c_3t^{1/3}(\alpha^{-1}+\tau^n \xi)+\tau^n\z}.
		 	\label{alphaproduct}
		 \end{equation}
		 Each factor is $1+o(1)$ since, for some $\epsilon>0$, for all $n\geq0$, $\alpha^{-1}+\tau^n\xi>\epsilon$ (this is where we need the assumption $\s<\rho^2$).
		 
		 Notice that near $z=1$,
		 \begin{equation*}
		 	f(\mu,z)=O\left(\frac{1}{|1-z|}\right) \text{ and } f(\mu,z)=\frac{\mu^{-1}}{1-z}+O(1),
		 \end{equation*}
		 so $\mu f(\mu,\z/\eta)/\eta$ in the kernel of $B_1$ becomes
		 \begin{equation}
		 	O\left(\frac{1}{|\eta-\zeta|}\right)\text{ and } \frac{1}{\eta-\z}+O(t^{-1/3})
		 	\label{fterm}
		 \end{equation} 
		 after the substitutions. From (\ref{cubic}), we have that, for some $\delta>0$, $e^{\psi(\z)}$ and $e^{-\psi(\eta)}$ are, respectively, $O(e^{-\delta|\z|^3})$ and $O(e^{-\delta|\eta|^3})$ after rescaling. Thus we can bound the rescaled kernels by constants times
		 \begin{equation*}
		 	\frac{e^{-\delta|\z|^3}}{|\z-\eta|}, \hspace{4mm}\frac{e^{-\delta|\eta|^3}}{|\eta-\z|},
		 	\label{bounds}
		 \end{equation*}
		 respectively. These are Hilbert-Schmidt operators, so to have convergence of the operators in Hilbert-Schmidt norm (and thus trace norm convergence of their product) it is enough to have pointwise convergence of their kernels. 
		 
		 The error in (\ref{fterm}) goes to zero and, since we may assume $a>1/4$, so does the error term in (\ref{cubic}). Thus the kernels have pointwise limits
		 \begin{equation*}
		  	\frac{e^{-\z^3/3+s\z}}{\z-\eta}, \hspace{4 mm} \frac{e^{\eta^3/3-s\eta}}{\eta-\z},
		 \end{equation*}
		 respectively. Thus, $T_1$ converges to the operator $J^{(0)}$ with kernel 
		 \begin{equation}
		 	J^{(0)}(\eta,\eta')=\int_{\G_\z} \frac{e^{-\z^3/3+s\z+(\eta')^3/3-s\eta'}} {(\z-\eta)(\eta'-\z)} d\z.
		 	\label{J0limit}
		 \end{equation}
		
		The contour $\G_\z$ converges to the rays from $-c_3$ to $-c_3+\infty e^{\pm 2\pi i/3}$ and the contour $\G_\eta$ converges to the rays from $0$ to $\infty e^{\pm \pi i/3}$.
		
		Since $\text{Re}(\z-\eta')<0$ when $\z\in\G_\z$ and $\eta'\in\G_\eta$, it follows that
		\begin{equation*}
			\frac{e^{s(\z-\eta')}}{\eta'-\z}=\int_s^\infty e^{x(\z-\eta')}dx.
		\end{equation*}
		Thus (\ref{J0limit}) is equal to
		\begin{equation*}
			\int_s^\infty\int_{\G_\z}\frac{e^{-\z^3/3+(\eta')^3/3+x(\z-\eta')}} {\z-\eta}d\z dx.
		\end{equation*}
		Thus, the limiting operator can be written as the product of an operator from $L^2(s,\infty)$ to $L^2(\G_\eta)$ with kernel
		\begin{equation*}
			\int_{\G_\z}\frac{e^{-\z^3/3+x\z}}{\z-\eta}d\z
		\end{equation*}
		and an operator from $L^2(\G_\eta)$ to $L^2(s,\infty)$ with kernel
		\begin{equation*}
			e^{-x\eta+\eta^3/3}.
		\end{equation*}
		Changing the order of the operators preserves the Fredholm determinant, and
		\begin{equation*}
			\int_{\G_\z}\int_{\G_\eta}\frac{e^{-\z^3/3+\eta^3/3+y\z-x\eta}} {\z-\eta}d\eta d\z 
			=-K_{\text{Airy}}(x,y)
		\end{equation*}
		 so we can conclude
		 \begin{equation}
		 	\det(I+\mu T_1)\to\det(I-K_\text{Airy}\chi_{(s,\infty)})=F_2(s).
		 	\label{T_1}
		 \end{equation}
		 
		 The operator $T_2$ can be written as
		 \begin{equation*}
		 	A_2B_2\sum_{j=1}^L\frac{w_j}{w_j\xi-\tau}c_3^{-1}t^{-1/3}
		 \end{equation*}
		 where, before substitution, $A_2:L^2(\G_\z)\to L^2(\G_\eta)$ and $B_2:L^2(\G_\eta)\to L^2(\G_\z)$ have
		 kernels
		 \begin{equation*}
		 \frac{e^{\psi(\z)}}{\z-\eta}\prod_{n=0}^\infty \frac{1+\tau^n\alpha\eta}{1+\tau^n\alpha\z}, \hspace{4 mm} (\z-\eta)\frac{\mu f(\mu,\z/\eta)}{\eta e^{\psi(\eta)}}
		 \end{equation*}
		 respectively. Notice $A_2=A_1$ and $B_2$ is bounded by $e^{-\delta|\eta|^3}$, which is Hilbert-Schmidt, so $A_2B_2$ converges in trace norm to its pointwise limit $J^{(1)}$, with kernel
		 \begin{equation*}
		 J^{(1)}(\eta,\eta')=-\int_{\G_\z} \frac{e^{-\z^3/3+s\z+(\eta')^3/3-s\eta'}} {\z-\eta} d\z.
		 \label{J1limit}
		 \end{equation*}
		 Notice that $J^{(1)}=\frac{d}{ds}J^{(0)}$.
		 
		 For the operator $T_3$, consider 
		 \begin{equation*}
		 	E(\z,\eta';w)=\sum\z^nh_k(\eta')
		 \end{equation*}
		  where $h_k(\eta')=O(t^{-2/3}|\eta'|^{k+1})$ and the (finite) sum ranges over appropriate values of $n,k$. Thus we have
		 \begin{equation*}
		 	T_3(\eta,\eta')=\sum\int_{\G_\z}\big(A_1(\eta,\z)\z^n\big)\big(B_1(\z,\eta')h_k(\eta')\big)d\z
		 \end{equation*}
		 where $A_1(\eta,\z)$, $B_1(\z,\eta')$ are as in (\ref{AB}). From the bounds given in (\ref{bounds}), we see that $A_1(\eta,\z)$ and $B_1(\z,\eta')$ are decaying exponentially in $\z,\eta'$, respectively, and so after multiplying the kernels by $|\z|^n$ and $|\eta'|^k$, respectively, they remain Hilbert-Schmidt. Thus the trace norm of $T_3$ is $O(t^{-2/3})$.
		 
		 We have shown that $T_1=J^{(0)}+o(1)$ (where the $o(1)$ term is independent of all $w_j$), $T_2=J^{(1)}\wsum c_3^{-1}t^{-1/3}+o(t^{-1/3})$, and $T_3=O(t^{-2/3})$. This gives us
		 \begin{flalign*}
		 \det(I+\mu J_{x,L,m,\rho}(w))
		 =\det\left(I+J^{(0)}+o(1)+J^{(1)}\sum_{j=1}^L\frac{w_j}{w_j\xi-\tau}c_3^{-1}t^{-1/3}+o(t^{-1/3})\right)\\
		 =\det(I+J^{(0)}+o(1))\det\left(I+(I+J^{(0)})^{-1}J^{(1)} \sum_{j=1}^L\frac{w_j}{w_j\xi-\tau}c_3^{-1}t^{-1/3}+o(t^{-1/3})\right)\\
		 =(F_2(s)+o(1))\left[1+\text{tr}((I+J^{(0)})^{-1}J^{(1)})\sum_{j=1}^L\frac{w_j}{w_j\xi-\tau}c_3^{-1}t^{-1/3}+o(t^{-1/3})\right]
		 \end{flalign*}
		 Recall $J^{(1)}=\frac{d}{ds}J^{(0)}$, so it follows that
		 \begin{equation*}
		 	\tr((I+J^{(0)})^{-1}J^{(1)})=\frac{d}{ds}\log\det(I+J^{(0)})=\frac{F_2'(s)}{F_2(s)},
		 \end{equation*}
		 therefore
		 \begin{equation}
		 \det(I+\mu J_{x,L,m,\rho}(w))=F_2(s)+o(1)+F_2'(s)\wsum c_3^{-1}t^{-1/3}+o(t^{-1/3}).
		 \label{Jdet}
		 \end{equation}
		 The $o(1)$ is independent of the $w_j$, as before. All terms are independent of $\mu$, so to evaluate the $\mu$-integral in (\ref{J-rep}), it suffices to compute
		 \begin{equation}
		 	\int(\tau^L\mu\tau)_\infty\frac{d\mu}{\mu^L}=(-1)^{L-1}\frac{\tau^{(L-1)(3L-2)/2}}{(1-\tau)\cdots(1-\tau^{L-1})}.
		 	\label{mu}
		 \end{equation}
		 From section IV of \cite{blocks}, we know
		 \begin{equation}
		 	\int_{\Gamma_{0,\tau}}\cdots\int_{\Gamma_{0,\tau}}\prod_{j=1}^{L}\frac{(w_j-1)^{L-j}}{w_j(w_j-\tau)^{L-j+1}}\prod_{i<j}\frac{w_j}{w_j-\tau w_j}dw_L\cdots dw_1=0
		 	\label{wintanalytic}
		 \end{equation}
		 and
		 \begin{equation*}
		 \int_{\Gamma_{0,\tau}}\cdots\int_{\Gamma_{0,\tau}}\prod_{j=1}^{L}\frac{(w_j-1)^{L-j}}{w_j(w_j-\tau)^{L-j+1}}\prod_{i<j}\frac{w_j}{w_j-\tau w_j}\wsum
		 dw_L\cdots dw_1
		 \end{equation*}
		 \begin{equation}
		 	=-\frac{\xi^{L-1}}{(1-\xi)^L}\frac{(1-\tau)\cdots(1-\tau)^L}{\tau^{L^2}}.
		 	\label{wintsum}
		 \end{equation}
		 This, together with (\ref{mu}) and (\ref{Jdet}), shows
		 \begin{equation*}
		 	\bbP_{L,\rho}(x,m,t)=(-1)^{L-1}c_3^{-1}\frac{\xi^{L-1}}{(1-\xi)^L}F_2'(s)t^{-1/3}+o(t^{-1/3}).
		 \end{equation*}
		 Since $\xi=-\sqrt{\s}/(1-\sqrt{\s})$, $c_2=\s^{-1/6}(1-\sqrt{\s})^{2/3}$ and $c_3=\s^{-1/6}(1-\sqrt{\s})^{5/3}$, we have that $c_3^{-1}/(1-\xi)=c_2^{-1}$ and $\xi/(1-\xi)=-\sqrt{\s}$, allowing us to conclude
		 \begin{equation*}
		 	\bbP_{L,\rho}(x,m,t/\g)=c_2^{-1}\s^{(L-1)/2}F_2'(s)t^{-1/3}+o(t^{-1/3})
		 \end{equation*}
		 when $\s<\rho^2$.
		 
		 \section{Asymptotics when $\s=\rho^2$} \label{case2}
		 
		 When $\s=\rho^2$, $\xi=-\alpha^{-1}$, so we set $\s=\rho^2+o(t^{-1/3})$, making $\xi=-\alpha^{-1}+o(t^{-1/3})$. To avoid the singularity at $-\alpha^{-1}$, we deform the $\eta$- and $\z$-contours used in the case where $\s<\rho^2$ so that instead of passing through $\xi$ and $\xi-t^{-1/3}$ they pass through $\xi+2\delta t^{-1/3}$ and $\xi+\delta t^{-1/3}$, respectively ($\delta>0$ is arbitrary and fixed). After making the substitutions (\ref{sub}) in a neighborhood of $\xi$, the limiting $\eta$-contour $\G_\eta$ consists of the rays from $2\delta c_3$ to $2\delta c_3+\infty e^{\pm \pi i/3}$ and the limiting $\z$-contour $\G_\z$ consists of the rays from $\delta c_3$ to $\delta c_3 +\infty e^{\pm2\pi i/3}$. The proof goes almost exactly the same as the $\s<\rho^2$ case, except that in (\ref{alphaproduct}), the factors for $n>0$ are still $1+o(1)$, but the $n=0$ factor is
		 \begin{equation*}
		 	\frac{c_3t^{1/3}(\alpha^{-1}+\xi)+\eta}{c_3t^{1/3}(\alpha^{-1}+\xi)+\z}=\frac{\eta}{\z}+o(1)
		 \end{equation*}
		 since $\alpha^{-1}+\xi=o(t^{-1/3})$. This changes the kernel of $J^{(0)}$ to
		 \begin{equation*}
		 J^{(0)}(\eta,\eta')=\int_{\G_\z} \frac{e^{-\z^3/3+s\z+(\eta')^3/3-s\eta'}} {(\z-\eta)(\eta'-\z)}\frac{\eta}{\z} d\z,
		 \end{equation*}
		 and the kernel of $J^{(1)}$ to
		 \begin{equation*}
		 J^{(1)}(\eta,\eta')=\int_{\G_\z} \frac{e^{-\z^3/3+s\z+(\eta')^3/3-s\eta'}} {(\z-\eta)}\frac{\eta}{\z} d\z.
		 \end{equation*}
		 By the same reasoning as in the $\s<\rho^2$ case, $J^{(0)}$ has the same Fredholm determinant as an operator with kernel
		 \begin{equation*}
		 \int_{\G_\z}\int_{\G_\eta}\frac{e^{-\z^3/3+\eta^3/3+y\z-x\eta}} {\z-\eta}\frac{\eta}{\z}d\eta d\z
		 \end{equation*}
		 \begin{equation*}
		 =\int_{\G_\z}\int_{\G_\eta}\frac{e^{-\z^3/3+\eta^3/3+y\z-x\eta}} {\z-\eta}d\eta d\z +\int_{\G_\eta}e^{\eta^3/3-x\eta}d\eta\int_{\G_\z}\frac{e^{-\z^3/3+y\z}} {\z} d\z
		 \end{equation*}
		 \begin{equation*}
		 	=-\left[K_{\text{Airy}}(x,y)+\Ai(x)\int_{-\infty}^y\Ai(z)dz\right].
		 \end{equation*}
		 Thus, the determinant of $J^{(0)}$ is $F_1(s)^2$. We still have $\frac{d}{ds}J^{(1)}=J^{(0)}$ and so we can use (\ref{mu}), (\ref{wintanalytic}), and (\ref{wintsum}) to conclude that, when $\s<\rho^2$,
		 \begin{equation*}
		 \bbP_{L,\rho}(x,m,t/\g)=c_2^{-1}\s^{(L-1)/2}\frac{d}{ds}(F_2(s)^2)^{-1/3}+o(t^{-1/3}).
		 \end{equation*}

	\end{appendices}
	
	\section*{Acknowledgments}
	The author thanks Harold Widom for providing the proof in Appendix \ref{appwint}, and Craig Tracy for helping in all stages of this work.
	This work was supported by the National Science Foundation through grants DMS-1207995 and DMS-1809311.

\end{document}